\documentclass[reqno, a4paper,oneside]{amsart}
\usepackage{amsthm,amsmath,amssymb, amscd}
\usepackage[]{geometry}
\usepackage[colorlinks=true,allcolors=blue!80!black]{hyperref}
\usepackage[english]{babel}
\usepackage{graphicx}
\usepackage{tikz}
\usepackage{tikz-3dplot,pgfplots}
\usetikzlibrary{decorations.markings} % for oriented edges in tikz
\usepackage{color}
\usepackage{enumerate}
\usepackage{mathtools}
\usepackage[color=pink, textwidth=25mm]{todonotes}
\usepackage{kbordermatrix}
\renewcommand{\kbldelim}{(}% Left delimiter
\renewcommand{\kbrdelim}{)}% Right delimiter

\usepackage{nicematrix}
\usepackage{float}

\theoremstyle{definition}
\newtheorem{definition}{Definition}[section]
\newtheorem{remark}[definition]{Remark} 
\newtheorem{defn}[definition]{Definition} 
\newtheorem{exam}[definition]{Example}

\newtheorem{question}{Question}
\newenvironment{q2prime}
  {
   \begin{question}}
  {\end{question}}

\theoremstyle{plain}

\newtheorem{thrm}[definition]{Theorem} 
\newtheorem{prop}[definition]{Proposition}
\newtheorem{lem}[definition]{Lemma}
\newtheorem{cor}[definition]{Corollary}

\newtheorem{abcthm}{Theorem}
\newtheorem{abccor}[abcthm]{Corollary}

\newtheorem*{conj}{Conjecture}

\newcommand{\Z}{\mathbb{Z}}
\newcommand{\T}{\mathbb{T}}
\newcommand{\G}{\Gamma}

\newcommand{\Sal}{\operatorname{Sal}}
\newcommand{\supp}{\operatorname{supp}}
\newcommand{\stc}{\operatorname{sc}}
\newcommand{\rk}{\operatorname{rank}}
\newcommand{\gen}{\operatorname{gen}}
\newcommand{\al}{\alpha}

\makeatletter
\newcommand{\extp}{\@ifnextchar^\@extp{\@extp^{\,}}}
\def\@extp^#1{\mathop{\bigwedge\nolimits^{\!#1}}}
\makeatother

\newcommand{\comment}[1]{}

\newcommand{\fakeenv}{} %%% prints the emptystring
\newenvironment{restate}[2]  %%% restate takes two arguments
{
  \renewcommand{\fakeenv}{#2} %%% So now \fakeenv prints #2
  \theoremstyle{plain}
  \newtheorem*{\fakeenv}{#1~\ref{#2}} %%% so now #2 is the name of a
  \begin{\fakeenv}  
}
{
  \end{\fakeenv}
}

\title{The minimal genus problem for right angled Artin groups}
\author[]{Rachael Boyd}
\address{University of Cambridge\\ DPMMS, Wilberforce Road, Cambridge CB3 0WB, UK}
\email{rachaelboyd@dpmms.cam.ac.uk}
\author[]{Thorben Kastenholz}
\address{University of G\"ottingen, Bunsenstrasse 3-5, 37073 G\"ottingen, Germany}
\email{thorben.kastenholz@mathematik.uni-goettingen.de}
\author[]{Jean Pierre Mutanguha}
\address{Max Planck Institute for Mathematics, Bonn}
\email{mutanguha@mpim-bonn.mpg.de}

\makeatletter
\@namedef{subjclassname@2020}{%
     \textup{2020} Mathematics Subject Classification}
\makeatother

\begin{document}


\begin{abstract}
    We investigate the minimal genus problem for the second homology of a right angled Artin group (RAAG). Firstly, we present a lower bound for the minimal genus of a second homology class, equal to half the rank of the corresponding cap product matrix. We show that for complete graphs, trees, and complete bipartite graphs, this bound is an equality, and furthermore in these cases the minimal genus can always be realised by a disjoint union of tori. Additionally, we give a full characterisation of classes that are representable by a single torus. %However, it is not true in general that the minimal genus of a second homology class of a RAAG is necessarily realised by a disjoint union of tori: we construct a genus two representative for a class in the pentagon RAAG.
    However, the minimal genus of a second homology class of a RAAG is not always realised by a disjoint union of tori as an example we construct in the pentagon shows.
\end{abstract}
\maketitle

\section{Introduction}
In this paper, we investigate the minimal genus of a second homology class of a right angled Artin group.
We always consider integral homology, and say that a continuous map~$f\colon \Sigma \to X$ from a compact (potentially disconnected) oriented surface \emph{represents a second homology class~$\alpha$} of a space $X$ if the induced map on homology sends the fundamental class~$[\Sigma] \in H_2(\Sigma)$ to~$\al$. %Recall that the fundamental class of~$\Sigma$ is the class $[\Sigma] \in H_2(\Sigma)$ determining the orientation. 
We define the \emph{minimal genus} of~$\alpha$, denoted by $\gen(\alpha)$, to be the minimal genus of a surface~$\Sigma$ representing~$\al$ in this way, where the genus of a disconnected surfaces is defined as the sum of the genera of the connected components. When we talk about the minimal genus of a second homology class of a group~$G$, we mean the minimal genus of a second homology class in~$BG$---the classifying space for the group. Note that homotopy equivalences preserve the minimal genus, hence any model for the classifying space yields the same minimal genus. We will restrict ourselves to the case where~$G$ is a right angled Artin group, also known as a RAAG.

Recall that a right angled Artin group~$A_\G$, is a group associated to a (finite) simple graph~$\G$ whose generators are given by the vertices~$V(\G)$ and commuting relations by the edges~$E(\G)$, i.e.~the presentation is:
\[ A_\G = \langle\,V(\G)~|~st = ts~\forall\,\{s,t\} \in E(\G)\,\rangle. \]
Extreme examples are the free abelian group~$\Z^n$, corresponding to $\G$ being the complete graph on $n$ vertices, and the free group~$F_n$, corresponding to $\G$ being $n$ disjoint vertices. 

The second integral homology $H_2(A_\G)$ can be identified with $\Z^{E(\G)}$, and we call the \emph{support} of the homology class~$\al \in H_2(A_\G)$ the union of edges whose corresponding entries in $\al \in \Z^{E(\G)}$ (as a vector) are not zero. In the specific case of RAAGs, the minimal genus of~$\al \in H_2(A_\G)$ is bounded above by the number of edges in the support of~$\al$.

There is a very general lower bound for the minimal genus of a second homology class. Namely, given a second homology class $\alpha \in H_2(X)$ of a topological space $X$, consider the \emph{cap product map}~$\alpha \cap - \colon H^1(X) \to H_1(X)$. It follows from naturality of the cap product that the image of~$\alpha \cap -$ must lie in the image of the first homology of any representative of $\alpha$. Since the genus is half the rank of the first homology of a surface, this yields the \emph{cap product inequality}:
\[
    {\frac{1}{2}} \, \rk(\alpha \cap -) \leq \gen(\alpha).
\]
The anti-symmetry of the cup product implies that $\rk(\alpha \cap -)$ is always an even integer. 

This inequality is in general far from an equality (as an example take the classifying space for any perfect group with non vanishing second homology), but we were able to show that it is indeed an equality for large families of RAAGs. 

\subsection{Results}
It was shown in \cite{KastenholzPedron} that the minimal genus only depends on the fundamental group, in the sense that the minimal genus of~$\al \in H_2(X)$ is the same as the minimal genus of the image of~$\al$ in~$H_2(\pi_1(X))$. All of our results could therefore be phrased in terms of second homology classes of spaces with fundamental group the specified RAAG. However since in practise we prove these statements for a classifying space, we state the results in terms of group homology.

For any RAAG, we introduce a diagrammatic description of a class $\al \in H_2(A_\G)$ and this provides a matrix description of $\al \cap -$, called the \emph{connection matrix} and denoted by~$M_\al$.

The case where $\G$ is a complete graph, i.e.~$A_\G \cong \mathbb{Z}^n$, serves as a guiding example for the minimal genus problem for all RAAGs. Since every separating curve in a surface is a commutator in~$\pi_1$, it follows that the minimal genus for any space with abelian fundamental group will be realised by a disjoint union of tori. 
Using this, we can translate the minimal genus problem for $\mathbb{Z}^n$ to an algebraic problem about skew-symmetric integer matrices. We obtain the following:

\begin{restate}{Theorem}{abcthm - torus}
Let $\G$ be a complete graph---i.e.~$A_\G\cong\Z^n$ and the~$n$-torus is a model for the classifying space---and $\al \in H_2(A_\G)$. 
Then the minimal genus $\gen(\al)$ is equal to the cap bound $\frac{1}{2} \, \rk(M_\alpha)$. Furthermore, the minimal genus is always realised by a disjoint union of tori.
\end{restate}

This complete solution to the minimal genus problem for $\Z^n$ leads to the following questions, which we (partially) answer in this paper:
\begin{question}\label{question - cap bound}
    Is the cap product inequality always an equality for a RAAG?
\end{question}
\begin{question}\label{question - realised by dj tori}
    Does every class in the second homology of a RAAG have a minimal genus representative that is a disjoint union of tori?
\end{question}

We were able to answer both questions for two large classes of RAAGs:
\begin{restate}{Theorem}{abcthm - cap bound equality}
Let $\G$ be a complete bipartite graph or a tree and $\al \in H_2(A_\G)$. Then the minimal genus $\gen(\al)$ is equal to the cap bound $\frac{1}{2} \rk(M_\alpha)$. Furthermore, the minimal genus can always be realised by a disjoint union of tori.
\end{restate}

By Proposition~\ref{prop - connected components cap bound}, if the cap product inequality is an equality in each of the connected components of a graph, then it is an equality for the whole graph. Hence the statement of Theorem~\ref{abcthm - cap bound equality} holds for disjoint unions of complete graphs, complete bipartite graphs, and trees.   

In the 1980s, Droms~\cite{Droms} classified all RAAGs that appear as fundamental groups of 3-manifolds: $\G$ must be a disjoint union of trees and triangles. Special cases of Theorems~\ref{abcthm - torus} and \ref{abcthm - cap bound equality} come together with Corollary~3.6 in \cite{KastenholzPedron} to give  us the following corollary.

\setcounter{abcthm}{2}
\begin{abccor}\label{abccor-droms}
    Let~$X$ be a 3-manifold such that~$\pi_1(X)$ is a RAAG. Then for any $\al \in H_2(X)$, the minimal genus $\gen(\al)$ is equal to the cap bound $\frac{1}{2} \rk(M_\alpha)$. Furthermore, the minimal genus can always be realised by a disjoint union of tori.
\end{abccor}
\setcounter{abcthm}{0}

Although we were not able to answer Question~\ref{question - cap bound} in general, we showed that the cap bound completely determines which classes are representable by a single torus.

\begin{restate}{Theorem}{abcthm - n partite torus}
Let~$\G$ be any simple graph. Then a nontrivial second homology class $\alpha \in H_2(A_\G)$ is representable by a torus if and only if $\rk(M_\alpha) = 2$. Furthermore, the support of such a homology class is a complete n-partite graph.
\end{restate}

These results provide a possible, albeit cumbersome, way to compute the minimal number of disjoint tori one needs to represent a second homology class: cover the support by complete~$n$-partite graphs. 

In all of the above cases, the minimal genus was always realised by a disjoint union of tori. However, we give a negative answer to Question~\ref{question - realised by dj tori} in general:

\begin{restate}{Theorem}{abcthm - counterexample to q2}
There exists a RAAG~$A_\G$ with a second homology class whose minimal genus representative cannot be realised by a disjoint union of tori.
\end{restate}

To prove this theorem, we find a second homology class, when~$\G$ is the pentagon, that has minimal genus two, but is not representable by fewer than three disjoint tori. 

\subsection{Future directions and context}
One obvious future direction would be to answer Question~\ref{question - cap bound} in full generality. Although we were unable to do this using the tools we developed, we conjecture the following (partly because our search for a counterexample was fruitless):

\begin{conj}
    The cap product inequality is an equality for any RAAG.
\end{conj}

In answering Question~\ref{question - realised by dj tori}, the representative we construct for Theorem~\ref{abcthm - counterexample to q2} is not $\pi_1$-injective. It is natural to then ask if the failure of $\pi_1$-injectivity is a necessary condition.

\begin{question}\label{question - surface subgroups}
    Does there exist a second homology class of a RAAG with a minimal genus representative that is not a disjoint union of tori but is injective on fundamental groups?
\end{question}

This question might be of interest to people studying surface subgroups in RAAGs. 
Fortunately, other examples of minimal genus representatives that fail to be $\pi_1$-injective like the example in Theorem~\ref{abcthm - counterexample to q2} could also be interesting in their own right. We remark that a minimal genus representative with the maximal number of connected components cannot map an essential simple closed curve of the surface to a null-homotopic loop: otherwise, performing surgery at this curve would either increase the number of components (if the curve is separating) or decreases the genus (if the curve is non-separating) while preserving the second homology class represented by the map. Thus for any minimal genus representative with the maximal number number of components, the induced map on fundamental groups is a homomorphism from a surface group to the RAAG that has no simple closed curves in its kernel.

Crisp, Sageev and Sapir asked whether any homomorphism of a surface group to a RAAG with no hyperbolic surface subgroups is necessarily injective if it has no simple closed curves in its kernel \cite[Problem~1.8]{Crisp}. The authors call the question an analogue of the ``simple curve in the kernel" problem for 3-manifolds. For instance, Stallings proved in~\cite{Stallings} that the restriction of the question to products of free groups is equivalent to the Poincar\'e Conjecture. To answer the general question, one might start by restricting Question~\ref{question - realised by dj tori} to RAAGs with no hyperbolic surface subgroups:

\begin{q2prime}
    In a RAAG with no hyperbolic surface subgroups, does every class in the second homology have a minimal genus representative that is a disjoint union of tori?
\end{q2prime}

Ideally, there would be a negative answer to this question that has the maximal number of components but is not $\pi_1$-injective---this would answer Crisp, Sageev and Sapir's problem.

Another future direction would be to investigate the minimal genus problem for other classes of Artin groups. By Proposition~3.5 in \cite{KastenholzPedron} and the fact that the second homotopy group of the Salvetti complex of any Artin group vanishes \cite[Proposition 1.13]{EliasWilliamson}, the minimal genus problem for an Artin group agrees with the minimal genus problem for the corresponding Salvetti complex. Unfortunately, there is no known formula for the second homology of a general Artin group. Akita and Liu \cite{AkitaLiu} give a general formula for the second homology with~$\Z/2\Z$ coefficients, but no integral results are known. Without such a result, a general investigation of the minimal genus problem for Artin groups seems to be out of reach.

\subsection{Outline}
In Section~\ref{section - preliminaries}, we provide some background on right angled Artin groups and the Salvetti complex. We also define the minimal genus and and introduce the descriptors we use to study it---the \emph{support} of a class and the corresponding \emph{connection matrix}. In Section~\ref{section - cap product inequality}, we construct our main tool---the \emph{cap bound inequality}---and prove some simple lemmas that we use throughout the paper, as well as Theorem~\ref{abcthm - torus}. Section~\ref{section - question 1} is devoted to answering Question~\ref{question - cap bound} for large classes of RAAGs---we prove Theorem~\ref{abcthm - cap bound equality} in this section. We then completely characterise which classes are representable by a single torus in Section~\ref{section - one torus}, proving Theorem~\ref{abcthm - n partite torus}. We finish in Section~\ref{section - question 2} with a negative answer to Question~\ref{question - realised by dj tori}, proving Theorem~\ref{abcthm - counterexample to q2}.

\subsection{Acknowledgements}
We would like to thank Mark Pedron for helpful conversations. The first and third authors also thank the Max Planck Institute for Mathematics in Bonn for its support and hospitality.

% \subsection{Data availability statement}
% Data sharing not applicable to this article.

\section{Preliminaries}\label{section - preliminaries}

\subsection{Right angled Artin groups}

We start by giving some background on right angled Artin groups. For a comprehensive introduction to RAAGs see the survey paper by Charney~\cite{Charney}.

\begin{defn}\label{defn- RAAG}
Every finite simple graph $\G$ with vertex set~$V(\G)$ and edge set $E(\G)$ determines a \emph{right angled Artin group}, or RAAG, $A_\G$, which is the group with presentation
\[
    A_\G =\langle V(\G) \mid st=ts \,\, \forall \{s,t\} \in E(\G) \rangle.
\] 
\end{defn}

\begin{exam}
Figure~\ref{fig - raag examples} shows a few examples of graphs $\G$ and their corresponding RAAGs.

\begin{figure}[ht]
    \begin{tikzpicture}[thick, scale=.84, decoration={
        markings,
        mark=at position 0.5 with {\arrow{>}}}]
		
		\path (0,0) -- node[below=.3, scale=1] {$A_\G\cong F_4$} (2,0);
		\path (0,0) -- node[left=.3, scale=1] {$\G = $} (0,2);
		
		\filldraw (0,0) circle [radius=0.1]
		          (2,2) circle [radius=0.1]
    		      (2,0) circle [radius=0.1]
		          (0,2) circle [radius=0.1];
		
		\begin{scope}[shift={(5,0)}]
		
    		\draw (0,0) -- node[below=.3, scale=1] {$A_\G\cong F_2 \times F_2$} (2,0);
    		\draw (0,0) -- node[left=.3, scale=1] {$\G = $} (0,2);
    		\draw (2,2) -- (2,0);
    		\draw (2,2) -- (0,2);
    		
    		\filldraw (0,0) circle [radius=0.1]
    		          (2,2) circle [radius=0.1]
    		          (2,0) circle [radius=0.1]
    		          (0,2) circle [radius=0.1];
    		          
		\end{scope}
		
		\begin{scope}[shift={(10,0)}]
		
    		\draw (0,0) -- node[below=.3, scale=1] {$A_\G\cong \Z^4$} (2,0);
    		\draw (0,0) -- node[left=.3, scale=1] {$\G = $} (0,2);
    		\draw (2,2) -- (2,0);
    		\draw (2,2) -- (0,2);
    		\draw (0,0) -- (2,2);
    		\filldraw[color=white] (1,1) circle [radius = 0.1];
    		\draw (0,2) -- (2,0);
    		
    		\filldraw (0,0) circle [radius=0.1]
    		          (2,2) circle [radius=0.1]
    		          (2,0) circle [radius=0.1]
    		          (0,2) circle [radius=0.1];
    		          
		\end{scope}
		
		\begin{scope}[shift={(3,-2)}]
		
    		\draw (0,0) node[left=.3, scale=1] {$\G = $} --
    		      (2,0) -- node[below=.3, scale=1] {$A_\G = \langle \, v_1, v_2, v_3, v_4~|~v_1v_2=v_2v_1, v_2v_3 = v_3v_2, v_3v_4=v_4v_3 \, \rangle$} 
    		      (4,0) -- (6,0);
    		
    		\filldraw (0,0) circle [radius=0.1] node[above, scale=1] {$v_1$} 
    		          (2,0) circle [radius=0.1] node[above, scale=1] {$v_2$} 
    		          (4,0) circle [radius=0.1] node[above, scale=1] {$v_3$} 
    		          (6,0) circle [radius=0.1] node[above, scale=1] {$v_4$};
    		          
		\end{scope}
	\end{tikzpicture}
    \caption{Four graphs on four vertices and their corresponding RAAGs.}
    \label{fig - raag examples}
\end{figure}
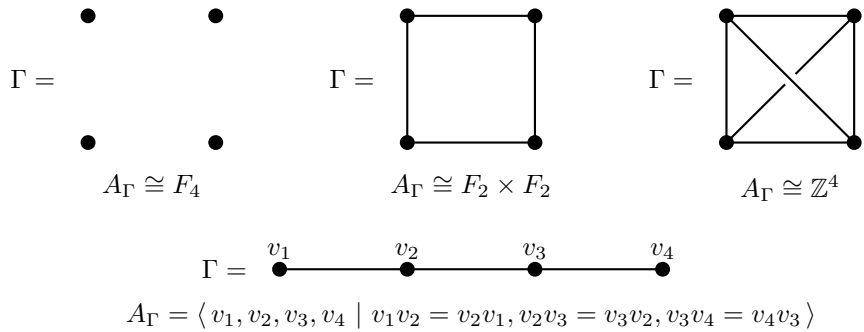
For an arbitrary simple graph~$\G$, the corresponding RAAG~$A_\G$ does not usually have a better description than the presentation given in Definition~\ref{defn- RAAG}.
\end{exam}

To study the minimal genus problem for RAAGs, we need two things: firstly a nice model for the classifying space of a RAAG, $BA_\G$, and secondly a method to describe a second homology class~$\al \in H_2(A_\G)$ since the minimal genus problem is concerned with such classes. For RAAGs, there exists a finite dimensional cube complex called the \emph{Salvetti complex} that is a model for the classifying space $BA_\G$. This complex was defined for general Artin groups by Salvetti in the 80s~\cite{Salvetti} and is a cube complex only for RAAGs. In general, it is not known whether it is always a classifying space, this is the well known \emph{$K(\pi,1)$ conjecture}. We therefore restrict to RAAGs for our definition. 

\begin{defn}\label{defn - salvetti}
Given a simple graph~$\G$ and corresponding RAAG~$A_\G$, the \emph{Salvetti complex} $\Sal_\G$ is the cube complex with:
\begin{itemize}
    \item one vertex, or 0-cube,~$x_0$;
    \item one edge, or $1$-cube, for each generator $s\in V(\G)$, attached to~$x_0$ at both ends;
    \item one square, or $2$-cube, for each edge $\{s_i,s_j\}\in E(\G)$, attached to the $1$-skeleton along the boundary edges using the relation $s_is_js_i^{-1}s_j^{-1}$. The image of each square is a $2$-torus in the $2$-skeleton;
    \item one $3$-cube, for each triangle in $\G$, attached to the $2$-skeleton by identifying opposite boundary squares with the $2$-tori corresponding to the three edges of the triangle; and
    \item generally, one $k$-cube, for each $k$-clique (complete graph on $k$ vertices) in $\G$, attached to the $(k-1)$-skeleton by identifying opposite boundary $(k-1)$-cubes with the $(k-1)$-tori corresponding to the $(k-1)$-cliques in the $k$-clique.
\end{itemize}
\end{defn}

\begin{exam}
In the extreme case when~$\G$ is a totally disconnected graph on $n$~vertices, the Salvetti complex~$\Sal_\G$ is a~\emph{rose}~$R_n= \bigvee_{i=1}^n S^1$ with $n$~petals. On the other hand, when~$\G$ is a complete graph on $n$~vertices, the Salvetti complex~$\Sal_\G$ is an $n$-torus $\T^n = (S^1)^n$. For an intermediate example, let~$\G$ be a square; then~$\Sal_\G$ is a product of roses~$R_2 \times R_2$. And as last example, suppose~$\G$ is a line with $3$~vertices; then~$\Sal_\G$ consists of two copies of a torus~$\T^2$ glued along a longitude in each copy.
\end{exam}

\begin{lem}[\cite{CharneyDavis}]\label{lem - salvetti is classifying space}
Let~$\G$ be a simple graph,~$A_\G$ the associated RAAG, and~$\Sal_\G$ the Salvetti complex. Then~$\Sal_\G$ is a model for the classifying space~$BA_\G$, or in other words, a~$K(A_\G,1)$.
\begin{proof}
The fundamental group of~$\Sal_\G$ is~$A_\G$ by construction. The fact that the universal cover is contractible follows from $\Sal_\G$ being a locally CAT(0) cube complex---a property unique to the Salvetti complex of a RAAG.
\end{proof}
\end{lem}

\begin{remark}[Orientation of 1-skeleton]
In Lemma~\ref{lem - salvetti is classifying space}, we implicitly choose an identification of~$\pi_1(\Sal_\G)$ with~$A_\G$. This identification automatically endows the Salvetti complex with a preferred orientation on its 1-skeleton. We choose such an identification now and fix it for the remainder of the paper.
\end{remark}

\begin{prop}
\label{prop - homology of salvetti}
Given a simple graph~$\G$, 
\begin{align*}
    H_1(A_\G) &= H_1(\Sal_\G) \cong \Z^{V(\G)} \quad \text{and}\\
    H_2(A_\G) &= H_2(\Sal_\G) \cong \Z^{E(\G)},
\end{align*}
where~$V(\G)$ and $E(\G)$ are the vertex set and edge set of~$\G$ respectively.
\end{prop}
This proposition follows immediately from the cellular chain complex of the Salvetti complex. But there is one small caveat, namely these isomorphisms implicitly choose orientations on the $1$-cells and $2$-cells of the Salvetti complex. We already chose an orientation on the $1$-cells in the above remark, and we will now address how to choose an orientation of the $2$-cells. We first orient $\G$:

\begin{defn}
Each edge $\{v,w\}\in E(\G)$ has two orientations given by the ordered pairs $(v,w)$ and $(w,v)$. An \emph{orientation} of $\G$ will be a choice of an oriented edge---either $(v,w)$ or $(w,v)$---for each edge $\{v, w\}$ in~$\G$. A simple graph with an orientation will be called an \emph{oriented graph}, and we denote the set of orientated edges by $E^{or}(\G)$. 
\end{defn}

\begin{remark}[Orientation of 2-skeleton]\label{rem - orientation of 2-skeleton}
Given an orientation of $\G$, we orient the $2$-cell in $\Sal_\G$ corresponding to an oriented edge $e=(v,w) \in E^{or}(\G)$ by considering the dual of $v^* \cup w^*$, where~$v^*$ and~$w^*\in H^1(\Sal_\G)$ denote the dual of the homology classes corresponding to $v$ and $w$ respectively (here we take the dual with respect to the basis given by the vertices in Proposition~\ref{prop - homology of salvetti}). Note that since~$v^* \cup w^*$ is a generator of~$H^2(\Sal_\G)$, its dual is an orientation of the corresponding torus, and thus gives an orientation of the 2-cell. Anti-symmetry of the cup product means that doing the same construction with the oriented edge~$(w,v)$ will yield the opposite orientation.
\end{remark}

Using an orientation on $\G$ and the induced orientation on $\Sal_\G$, we obtain a canonical generator for each~$\Z$ factor in~$H_2(A_\G)\cong \Z^{E(\G)}$---we let~$e_{(v,w)} \in \Z^{E(\G)}$ be the basis element corresponding to the oriented~$2$-cell given by the oriented edge~$(v,w) \in E^{or}(\G)$, and~$-e_{(v,w)}$ correspond to the same $2$-cell with the opposite orientation. This allows us to do the following.

Let~$\G$ be a simple graph and fix~$\al \in H_2(A_\G)$. To describe the class with a pictorial approach, we choose an orientation of the graph; this is equivalent to decorating each edge with an arrow. As mentioned in the preceding paragraph, an orientation of~$\G$ determines a basis for $H_2(A_\G) \cong \Z^{E(\G)}$.

\begin{defn}\label{defn-support} Let~$\alpha \in H_2(A_\G)$ be an arbitrary homology class, and~$l(v,w)$ be the integer coefficient for the basis vector~$e_{(v,w)}$ in the vector $\al \in \Z^{E(\G)}$.
The \emph{support} of~$\alpha$, denoted by~$\supp(\al)$, is the following labelled oriented graph: consider the oriented subgraph of~$\G$ spanned by oriented edges~$(v,w)$ where the label~$l(v,w)$ is non-zero and label the edges with the non-zero integer label~$l(v,w)$. 
\end{defn}

The support uniquely determines the class~$\alpha$ up to the following relation on labelled oriented graphs: 
\begin{center}
    \begin{tikzpicture}[thick, scale=.84, decoration={
        markings,
        mark=at position 0.5 with {\arrow{>}}}]
        
		\draw[postaction={decorate}] (0,0) -- node[above=2pt] {$n$} (2,0);
		\draw[postaction={decorate}] (6,0) -- node[above=2pt] {$-n$} (4,0);
		\draw (3,0) node {$=$};
			            
		\filldraw (0,0) circle [radius=0.1] node[above=1pt] {$v$}
		          (2,0) circle [radius=0.1] node[above=1pt] {$w$}
		          (4,0) circle [radius=0.1] node[above=1pt] {$v$}
		          (6,0) circle [radius=0.1] node[above=1pt] {$w$};
		
	\end{tikzpicture}
\end{center}
We omit the integer label on an edge when the label is zero. 
Note that the underlying graph of~$\supp(\al)$ does not depend on the chosen orientation of~$\G$ and we will sometimes also consider~$\supp(\al)$ as a subgraph of an unoriented graph~$\G$.

One of our main tools in this work is the following matrix, derived from the labelled support of a class $\al \in H_2(A_\G)$.

\begin{defn}\label{defn-connection matrix}
Let~$\G$ be a simple graph and~$\al \in H_2(A_\G)$ a given class. Following the preceding discussion, the class~$\al$ can be described by choosing an orientation of the graph~$\G$ and labelling each oriented edge $(v,w) \in E^{or}(\G)$ with an appropriate integer $l(v,w)$.
The \emph{connection matrix} of the class~$\al$ is a square matrix~$M_\al$ with rows and columns indexed by the vertices of~$\G$, and whose matrix entries are given by: 
\[ (M_\al)_{v,w} = 
\begin{cases} 
0 & \text{if } v=w \\ 
l(v,w) & \text{if } (v,w) \in E^{or}(\G) \\
-l(w,v) & \text{if } (w,v) \in E^{or}(\G).
\end{cases} \]
\end{defn}

Note that the connection matrix is a skew-symmetric integer matrix and, due to the relation on labelled graphs preceding Definition~\ref{defn-support}, it is independent of the orientation of~$\G$ we choose when depicting~$\al$ as a labelled oriented graph.

We now give an example of a homology class, its labelled support, and its connection matrix.

\begin{exam} Let $\G$ be the square with vertices $\{ v_1, v_2, w_1, w_2 \}$ and an orientation given by $E^{or}(\G)=\{ v_1, v_2 \} \times \{ w_1, w_2 \}$. Then $H_2(A_\G)$ has basis given by $e_{(v_1, w_1)},\,e_{(v_1, w_2)},\,e_{(v_2, w_1)},$ and~$e_{(v_2, w_2)}$.

Set $\alpha = e_{(v_1, w_1)} - e_{(v_2, w_2)}$ and $\beta = 2 e_{(v_1, w_1)} + 4e_{(v_1, w_2)} + 3e_{(v_2, w_1)} + 6e_{(v_2, w_2)}$ in $H_2(A_\G)$. Then the connection matrices are
\[
M_\al =
\kbordermatrix{
 & v_1 & v_2 & w_1 & w_2 \\
v_1 & 0 & 0 & 1 & 0 \\
v_2 & 0 & 0 & 0 & -1  \\
w_1 & -1 & 0 & 0 & 0  \\
w_2 & 0 & 1 & 0 & 0} \quad \text{and} \quad 
M_\beta =
\kbordermatrix{
 & v_1 & v_2 & w_1 & w_2 \\
v_1 & 0 & 0 & 2 & 4 \\
v_2 & 0 & 0 & 3 & 6  \\
w_1 & -2 & -3 & 0 & 0  \\
w_2 & -4 & -6 & 0 & 0}.
\]
Figure~\ref{fig - a class} is a visual representation of the classes with their supports highlighted.

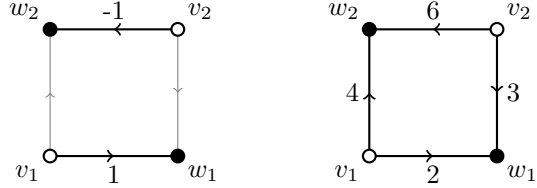
\begin{figure}[ht]
    \begin{tikzpicture}[thick, scale=.84, decoration={
        markings,
        mark=at position 0.5 with {\arrow{>}}}]
			
		\draw[postaction={decorate}] (0,0) -- node[below] {1} (2,0);
		\draw[color=gray, thin, postaction={decorate}] (0,0) -- (0,2);
		\draw[color=gray, thin, postaction={decorate}] (2,2) -- (2,0);
		\draw[postaction={decorate}] (2,2) -- node[above] {-1} (0,2);
		
		\draw[fill=white] (0,0) circle [radius=0.1] node[below left] {$v_1$}
		                  (2,2) circle [radius=0.1] node[above right] {$v_2$};
		\filldraw (2,0) circle [radius=0.1] node[below right] {$w_1$}
		          (0,2) circle [radius=0.1] node[above left] {$w_2$};
		
		\begin{scope}[shift={(5,0)}]
		
    		\draw[postaction={decorate}] (0,0) -- node[below] {2} (2,0);
    		\draw[postaction={decorate}] (0,0) -- node[left] {4} (0,2);
    		\draw[postaction={decorate}] (2,2) -- node[right] {3} (2,0);
    		\draw[postaction={decorate}] (2,2) -- node[above] {6} (0,2);
    		
    		\draw[fill=white] (0,0) circle [radius=0.1] node[below left] {$v_1$}
    		                  (2,2) circle [radius=0.1] node[above right] {$v_2$};
    		\filldraw (2,0) circle [radius=0.1] node[below right] {$w_1$}
    		          (0,2) circle [radius=0.1] node[above left] {$w_2$};
    		          
		\end{scope}
	\end{tikzpicture}
    \caption{An illustration of $\alpha$ and $\beta$.}
    \label{fig - a class}
\end{figure}
\end{exam}

\subsection{Minimal genus}
We will now define the minimal genus and discuss some its properties, in particular in the setting of group homology. See \cite{KastenholzPedron} for an introduction to the minimal genus problem.

\begin{defn}\label{defn-minimal genus}
Given a space~$X$ and a class~$\al \in H_2(X)$, we define the \emph{minimal genus} of~$\alpha$, denoted by $\gen(\alpha)$, to be the minimal genus of an compact oriented surface~$\Sigma$ such that there exists a continuous map~$f\colon \Sigma \to X$ and $f_*([\Sigma])=\al$ in~$H_2(X)$, where $f_*$ denotes the induced map on homology and $[\Sigma]$ is the fundamental class, i.e.~the orientation of~$\Sigma$. 
Here $\Sigma$ may have more than one connected component, and the surface genus of disconnected surfaces is the sum of the genera of the connected components. We say~$\al$ is \emph{representable by~$\Sigma$}.
\end{defn}

Since the disjoint union of two representatives represents the sum of the corresponding homology classes, we obtain the following lemma, which will be used later on when we split homology classes:
\begin{lem}\label{lem-subadditive}
    The minimal genus is subadditive:
    \[ \gen(\al + \beta) \leq \gen(\al)+\gen(\beta) \quad \forall \al, \beta \in H_2(X). \]
\end{lem}

In this work, we are interested in the minimal genus of a second integral homology class of a group---recall that the homology of a group~$G$ is defined to be the homology of its \emph{classifying space}~$BG$ \cite[Section 2.4]{Brown}. The classifying space is well-defined up to homotopy: it is a~$K(G,1)$ space. Homotopy equivalences preserve the minimal genus, hence any model for the classifying space yields the same minimal genus. 

\begin{defn}
The \emph{minimal genus} of a second homology class of a group $G$ is the minimal genus of a second homology class in~$BG$, the classifying space for the group. 
\end{defn}

In the case of RAAGs, we can restrict ourselves to studying the minimal genus of classes in~$H_2(\Sal_\G)$. More precisely, we define our classes via labelled oriented graphs and this concretely refers to a class in~$H_2(\Sal_\G)=\Z^{E(\G)}$, which has a canonical choice of basis as described in Remark~\ref{rem - orientation of 2-skeleton}. 

Moreover since the RAAGs corresponding to disjoint unions of graphs have Salvetti complexes given by wedges of the Salvetti complexes for the connected subgraphs, we will also need
the following proposition, which follows from \cite[Propistion~3.8]{KastenholzPedron}.
\begin{prop}
\label{proposition - genus wedge additive}
    Let $\G = \G_1 \sqcup \G_2$ denote a graph which is a disjoint union of two graphs. Then
    \begin{itemize}
    \item
        $H_2(\Sal_\G) \cong H_2(\Sal_{\G_1})\oplus H_2(\Sal_{\G_2})$
    \item
        For a class $H_2(\Sal_\G)\ni \alpha = \alpha_1 \oplus \alpha_2$ with $\alpha_i \in H_2(\Sal_{\G_i})$, we have 
        \[
            \gen(\alpha) = \gen(\alpha_1) + \gen(\alpha_2).
        \]  
    \end{itemize}
\end{prop}

\section{The cap product inequality}\label{section - cap product inequality}
In the first half of this section, we use the cap product to give a lower bound on the minimal genus. The second half uses this lower bound to compute the minimal genus of a second homology class in an $n$-torus~$\T^n$.

Given a space~$X$ and~$\alpha\in H_2(X)$, recall \cite[Section 3.3]{Hatcher} that the \emph{cap product} map of~$\alpha$ is the map
\[
\alpha \cap - \colon H^1(X)\to H_1(X).
\]

This map leads to a lower bound for the minimal genus:

\begin{prop}\label{prop-cap product inequality}
 Let~$X$ be any space and~$\alpha\in H_2(X)$. Then the following inequality holds
\[
2 \, \gen(\alpha)\geq \rk(\alpha \cap -).
\]
\begin{proof}
We first make a general observation. Suppose~$f:Y \to X$ is a map from some space $Y$ to~$X$ and let~$\beta \in H_m(Y)$ and~$\sigma \in H^n(X)$ be arbitrary classes. If~$f_*$ and $f^*$ are the induced maps on homology and cohomology respectively, then the cap product 
\[
f_*(\beta) \cap \sigma = f_*(\beta \cap f^*(\sigma))
\] 
lies in $f_*(H_{m-n}(Y))$. 

Now assume $\alpha \in H_2(X)$ and $f:\Sigma \to X$ represents $\alpha$, i.e.~$\Sigma$ is a possibly disconnected oriented surface of genus $g$ and~$f_*([\Sigma]) = \alpha$ where $[\Sigma] \in H_2(\Sigma)$ is the fundamental class.
Then, by the previous observation, the cap product $\alpha \cap \sigma$ lies in $f_*(H_1(\Sigma))$ for all $\sigma \in H^1(X)$. Therefore, the image of the cap product map $\alpha \cap -$ is contained in the image of $f_*:H_1(\Sigma) \to H_1(X)$, which has rank at most $2g$. 
So we get $2g \geq  \rk(\alpha \cap -)$. 
\end{proof}
\end{prop}

We now consider the case where~$\G$ is a simple graph and~$\alpha\in H_2(A_\G)$. In this case, it is easy to compute the bound~$\rk(\alpha \cap -)$ using the connection matrix~$M_\alpha$ (introduced in Definition~\ref{defn-connection matrix}).

\begin{prop}\label{prop-cap bound matrix rank}
Let $\G$ be a simple graph,~$\alpha\in H_2(A_\G)$ an arbitrary class, and~$M_\alpha$ its connection matrix. 
Then the connection matrix~$M_\alpha$ is also the matrix representation of the cap product map
\[
\alpha \cap - \colon H^1(A_\G) \to H_1(A_\G)
\]
with respect to the basis given by the fixed orientation on the $1$-skeleton.
\begin{proof} 
    The $2$-skeleton of the Salvetti complex $\Sal_{\Gamma}$ is a quotient of the space 
    \[X = \bigsqcup_{\{v,w\} \in E(\Gamma)} (S^1\times S^1)_{\{v,w\}},\] where $(S^1\times S^1)_{\{v,w\}}$ is (a copy of) a $2$-torus. The quotient map $\pi \colon X \to \Sal_{\Gamma}^{(2)}$ and the inclusion $\Sal_\Gamma^{(2)}\to \Sal_\Gamma$ both induce isomorphisms on second homology. Furthermore, the inclusion $\Sal_\Gamma^{(2)}\to \Sal_\Gamma$ induces an isomorphism on homology and cohomology in degrees $1$ and $2$. A distinguished basis for the first homology of $\Sal_\Gamma^{(2)}$ is given by~$V(\Gamma)$, and a distinguished basis for the first homology of $X$ is given by $(S^1\times \{\ast\})_{\{v,w\}}$ and \emph{dual curves} $(\{\ast\} \times S^1)_{\{v,w\}}$ for each $\{v,w\} \in E(\Gamma)$. Additionally, a basis for the first cohomology of $\Sal_\Gamma^{(2)}$ is given by the duals~$v^*$ of every element $v$ of~$V(\Gamma)$. The quotient map~$\pi$ induces the following map on cohomology
    \begin{eqnarray*}
        \pi^* 
        \colon 
        H^1(\Sal_\Gamma^{(2)}) 
        &\to& 
        H^1(X)\\
        v^* &\mapsto& \sum_{p\in P_v} p^*
    \end{eqnarray*}
    Here~$P_v$ is the set of circles in~$X$ which map to the circle in~$\Sal_{\Gamma}$ corresponding to~$v$ under~$\pi$, thus each~$p\in P_v$ represents a class in~$H_1(X)$ and~$p^*$ denotes its dual in~$H^1(X)$.
   The cap product of a generator of $H_2((S^1\times S^1)_{\{s,t\}})$ in $H_2(X)$ with $\pi^*(v^*)$ is either the dual curve to the corresponding $p \in P_v$ or zero if there is no $p \in P_v$ that lies in that torus. In other words, if $(S^1\times S^1)_{\{s,t\}}$ is oriented and $p$ is $(S^1\times \{\ast\})_{\{s,t\}}$ then the cap product is given by
   \[
   [(S^1\times S^1)_{\{s,t\}}]\cap p^*=\pm [(\{\ast\}\times S^1)_{\{s,t\}}],
   \]where the sign depends on the orientation of the torus. Since the cap product is natural and linear, this yields the desired matrix~$M_\al$.
\end{proof}
\end{prop}
\begin{remark}
    For another proof to Proposition~\ref{prop-cap bound matrix rank}, consider the inclusion $\Sal_\Gamma \to (S^1)^{|V(\Gamma)|}$ coming from the abelianisation map. This map induces an isomorphism on first homology and an injection on second homology. We can deduce the result from the cap product structure of the torus. 
\end{remark}

Putting Proposition~\ref{prop-cap product inequality} and Proposition~\ref{prop-cap bound matrix rank} together, we get the \emph{cap product inequality:}

\begin{equation}\label{eq - cap bound}
\gen(\alpha)\geq \frac{1}{2} \, \rk(M_\alpha) \quad \text{for all } \alpha \in H_2(A_\G).
\end{equation}
The right-hand side,~$\frac{1}{2} \, \rk(M_\alpha)$, will be referred to as the \emph{cap bound}. We are interested in when the cap product inequality is an equality.

\begin{prop}\label{prop - connected components cap bound}
Let $\G$ denote a graph which is a disjoint union of graphs $\G_i$ and suppose that for all  $\G_i$ and all $\beta \in H_2(\Sal_{\G_i})$ we have $\gen(\beta) = \frac{1}{2}\rk(M_\beta)$, then $\gen(\alpha) = \frac{1}{2}\rk(M_\alpha)$ holds for all $\alpha \in H_2(\Sal_{\G})$.
\begin{proof}
    Note that if a graph has multiple connected components, then its connection matrix will be a block matrix. Hence the cap bound of a class $\alpha = \sum_i \alpha_i$, where the $\alpha_i$ come from projecting to the second homology of the connected components, will be the sum of the cap bounds of the $\alpha_i$. Combining this with Proposition~\ref{proposition - genus wedge additive} gives the result.
\end{proof}
\end{prop}

To conclude the section, we compute the minimal genus when~$A_\G\cong \Z^n$ using the cap bound. In this case, we have extra tools that we can use.  Firstly, there is an isomorphism of graded rings
\[H_*(\Z^n)\cong \bigwedge\nolimits^* \Z^n,
\]
where the ring structure on the left stems from the group multiplication in $\T^n$, i.e.~it is the Pontryagin ring structure. This means we can interpret $\al \in H_2(\Z^n)$ as a skew-symmetric bilinear form on~$\Z^n$.
Furthermore, taking the dual on the left, we have an isomorphism of graded rings
\[H^*(\Z^n)\cong \bigwedge\nolimits^* \Z^n,
\]
where the ring structure on the left is given by the cup product.
Under these two ring isomorphisms, the matrix representation~$M_\al$ of the cap product map
\[\alpha \cap - \colon \mathbb{Z}^n \cong H^1(\T^n) \to H_1(\T^n) \cong \mathbb{Z}^n\] yields the same matrix as interpreting $\alpha$ as a skew-symmetric bilinear form.

\begin{prop}\label{prop - n torus cap equality}
Let $\G$ be a complete graph on $n$ vertices and $\al \in H_2(A_\G) \cong \bigwedge^2 \Z^n$. 
Then the minimal genus $\gen(\al)$ is always realized by a disjoint union of tori. Furthermore, it is equal to the minimal number of elementary wedges, i.e.~elements of the form~$a \wedge b$, needed to represent~$\al$.
\begin{proof}
   A disjoint union of $2$-tori can serve as the minimal genus representative for any second homology class of the $n$-torus~$\T^n = (S^1)^n$ since the fundamental group $\pi_1(\T^n) \cong \Z^n$ is abelian. 
   
    Now suppose that we have a class in~$H_2(\T^n)$ that is representable by a torus, i.e.~by a map~$\tau \colon \T^2 \to \T^n$. Since the $2$-torus and the $n$-torus are aspherical, the map $\tau$ is, up to homotopy, determined by the induced homomorphism on the fundamental groups. Thus we may assume it is given by $(s,t) \mapsto g(s)\cdot h(t)$ for some pair of based loops $g,h \colon S^1\to \T^n$, where multiplication is the group operation in~$\T^n$. 
    Let $[g]$ and $[h]$ denote the images of~$g_*([S^1])$ and $h_*([S^1])$ respectively under the chain of isomorphisms $\pi_1(\T^n) \cong H_1(\T^n) \cong \mathbb{Z}^n$. Then under the identification 
    \[ H_2(\T^n) \cong  \extp^2 H_1(\T^n) \cong \extp^2 \mathbb{Z}^n,\] such a map $\tau$ sends a generator of $H_2(\T^2)$ to the elementary wedge $[g] \wedge [h]$. 
    
    Conversely, let an elementary wedge $a \wedge b \in \extp^2 \mathbb{Z}^n \cong H_2(\T^n)$ be given. Choose explicit representatives $\hat{a},\hat{b} \colon S^1 \to \T^n$ for $a,b\in \mathbb{Z}^n \cong H_1(\T^n) \cong \pi_1(\T^n)$. Then the map $\T^2 \to \T^n$ defined by $(s,t) \mapsto \hat{a}(s)\cdot \hat{b}(t)$ maps a generator of~$H_2(\T^2)$ to~$a \wedge b$, i.e.~$a \wedge b$ is representable by a torus.
    
    All in all, this proves that the minimal genus of any~$\alpha \in H_2(\Z^n) \cong \extp^2 \mathbb{Z}^n$ is the minimal number of elementary wedges, i.e.~elements of the form $a\wedge b$, needed to represent~$\alpha$.
\end{proof}
\end{prop}

Fortunately, we can use linear algebra (over the integers) to compute the minimal number of elementary wedges of an element in~$\extp^2 \Z^n$ using an appropriate skew-symmetric integer matrix.

\begin{prop}\label{prop - minimal number of elementary wedges}
Let $\G$ be a complete graph on $n$ vertices and $\al \in H_2(A_\G) \cong \extp^2 \Z^n$. 
Then the cap bound $\frac{1}{2} \, \rk(M_\alpha)$ equals the minimal number of elementary wedges needed to represent~$\al$.
\begin{proof}
    The solution to this purely algebraic problem is classical. By silently identifying $\mathbb{Z}^n$ with its dual, we can think of $\alpha \in  \extp^2 \mathbb{Z}^n$ as a skew-symmetric bilinear form on $\mathbb{Z}^n$ (see the discussion before Proposition~\ref{prop - n torus cap equality}). Following this, by Section~14 (Skew-symmetric Matrices) in \cite{FranklinVeblenIntegerMatrices}, such a skew-symmetric bilinear form---represented as a matrix~$M$---has a normal form~$\widehat{M}$ consisting of a direct sum of integer multiples of the standard hyperbolic form and the zero form. This normal form is shown below, where~$0 \neq \lambda_i\in \Z$ for~$1\leq i\leq k$.
\[
\widehat{M}=
\begin{pNiceArray}{ccccccc:ccc}[first-row, first-col]
&v_1 & w_1 & v_2 & w_2& \cdots&v_k&w_k& \Block{1-3}{\cdots}&\\
v_1 & 0 & \lambda_1 &&&\Block{4-3}<\huge>{0}&&&\Block{7-3}<\huge>{0}& \\
w_1& -\lambda_1 &0 &&&&&&&  \\
v_2&  && 0&\lambda_2&&  & &&\\
w_2&  && -\lambda_2&0 & &  & &&\\
\vdots&\Block{3-4}<\huge>{0}& && &\ddots& &&&\\
v_k && && && 0 & \lambda_k &&\\
w_k && && && -\lambda_k&0 &&\\
\hdottedline
\Block{3-1}{\vdots}& \Block{3-7}<\huge>{0}&& & &&&&\Block{3-3}<\huge>{0}&&\\
&&&&&&&&&\\
&&&&&&&&&
\end{pNiceArray}
\]
    
    Consider the basis of $\mathbb{Z}^n$ such that as a matrix~$\al$ has the above normal form (with~$k$ hyperbolic blocks), and let $v_1, w_1, \ldots, v_k, w_k$ denote the first~$2k$ basis vectors as shown in the matrix diagram, i.e.~paired according to the hyperbolic forms. Then it follows that $\alpha = \sum_{i=1}^k \lambda_i v_i \wedge w_i$, i.e.~it is a sum of~$k$ elementary wedges.
    By Proposition~\ref{prop - n torus cap equality}, we get $\gen(\al) \le k$.
    Evidently, the number of hyperbolic blocks $k$ equals half the rank of the matrix~$M$ representing the skew-symmetric form.
    Finally, the matrix representation~$M_\al$ of the cap-product $\alpha \cap - \colon \mathbb{Z}^n \cong H^1(\T^n) \to H_1(\T^n) \cong \mathbb{Z}^n$ yields the same matrix~$M$ as interpreting $\alpha$ as a skew-symmetric bilinear form. So, combining with the cap product inequality (Equation~\ref{eq - cap bound}), we get~$\gen(\al) = \frac{1}{2} \, \rk(M_\al)$.
\end{proof}
\end{prop}

Together, the previous two propositions give us our first major result:
\begin{abcthm}\label{abcthm - torus}
Let $\G$ be a complete graph---i.e.~$A_\G\cong\Z^n$ and the~$n$-torus is a model for the classifying space---and $\al \in H_2(A_\G)$. 
Then the minimal genus $\gen(\al)$ is equal to the cap bound $\frac{1}{2} \, \rk(M_\alpha)$. Furthermore, the minimal genus is always realised by a disjoint union of tori.
\end{abcthm}

\section{A partial answer to Question~\ref{question - cap bound}}\label{section - question 1}
The cap product inequality (Equation~\ref{eq - cap bound}) and Theorem~\ref{abcthm - torus} naturally lead us to ask:
\setcounter{question}{0}
\begin{question}
    Is the cap product inequality always an equality for a RAAG?
\end{question}
This section answers the question in the affirmative for two large classes of RAAGs:
\begin{abcthm}\label{abcthm - cap bound equality}
Let $\G$ be a complete bipartite graph or a tree and $\al \in H_2(A_\G)$. Then the minimal genus $\gen(\al)$ is equal to the cap bound $\frac{1}{2} \rk(M_\alpha)$. Furthermore, the minimal genus can always be realised by a disjoint union of tori.
\end{abcthm}

Since the proofs are quite different for a complete bipartite graph and a tree, we handle them in two separate subsections, culminating in Theorems~\ref{thrm - inequality is equality complete bipartite} and~\ref{thrm - inequality is equality trees} respectively.

\subsection{Bipartite graphs}\label{subsect - complete bipartite}
Unless otherwise stated, we assume in this subsection that $\G$ is a complete bipartite finite graph, i.e.~there is a partition of the vertices $V(\G) = \{ v_1, \ldots, v_n \} \sqcup \{ w_1, \ldots, w_m \}$ such that each edge in~$\G$ has endpoints $v_i$ and $w_j$ for some $i,j$ and, conversely, each pair $\{v_i, w_j\}$ corresponds to a (unique) edge in~$\G$.
Then $A_\G\cong F_n \times F_m$ where $F_n$ and $F_m$ are free groups generated by $\{ v_1, \ldots, v_n \}$ and $\{ w_1, \ldots, w_m \}$ respectively. 
Elements of $A_\G$ are considered as words $v w$ where $v \in F_n$ and $w \in F_m$ rather than as ordered pairs $(v, w)$. 
For an orientation on $\G$, we shall choose the oriented edges $E^{or}(\G) = \{v_1, \ldots, v_n \} \times \{w_1, \ldots, w_m \}$. We write $\overline{v_i}$ and  $\overline{w_j}$ for the images of~$v_i$ and $w_j$ in $H_1(F_n)$ and $H_1(F_m)$ respectively. In this setting, $H_2(A_\G) \cong H_1(F_n) \otimes H_1(F_m)$ is generated by $e_{(v_i, w_j)} = \overline{v_i} \otimes \overline{w_j}$  for $(v_i, w_j) \in E^{or}(\G)$. 

We start by constructing examples of classes representable by a torus.

\begin{lem}\label{Thorben}
Let~$\G$ be a complete bipartite graph and identify $A_\G \cong F_n \times F_m$. For any $v \in F_n$ and $w \in F_m$, the class $\al = \overline v \otimes \overline w \in H_2(A_\G)$ is representable by a torus.

\begin{proof}
Let~$\Z^2$ be generated by~$a$ and~$b$. For any $v \in F_n$ and $w \in F_m$, we define a homomorphism 
\begin{eqnarray*}
\tau:\Z^2 &\to& A_\G\\
a &\mapsto& v \\ 
b &\mapsto& w.
\end{eqnarray*}
Using the identifications $H_2(\Z^2) \cong \extp^2 \Z^2$ and $H_2(A_\G) \cong H_1(F_n) \otimes H_1(F_m)$, direct computation shows that the induced homomorphism $\tau_* \colon H_2(\Z^2) \to H_2(A_\G)$ maps the generator $a \wedge b$ to $\overline v \otimes \overline w$. As~$\T^2$ and~$\Sal_\G$ are $K(\pi,1)$-spaces, we get a map $\T^2 \to \Sal_\G$ that induces $\tau$ on fundamental groups, and hence the class $\al=\overline v \otimes \overline w$ is representable by a torus.
\end{proof}
\end{lem}

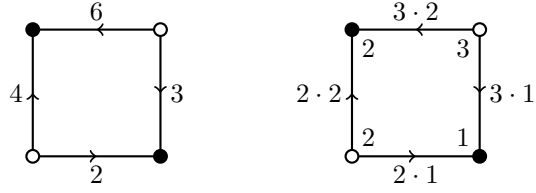
\begin{figure}[ht]
    \begin{tikzpicture}[thick, scale=.84, decoration={
        markings,
        mark=at position 0.5 with {\arrow{>}}}]
			
		\draw[postaction={decorate}] (0,0) -- node[below] {2} (2,0);
		\draw[postaction={decorate}] (0,0) -- node[left] {4} (0,2);
		\draw[postaction={decorate}] (2,2) -- node[right] {3} (2,0);
		\draw[postaction={decorate}] (2,2) -- node[above] {6} (0,2);
		
		\draw[fill=white] (0,0) circle [radius=0.1]
		                  (2,2) circle [radius=0.1];
		\filldraw (2,0) circle [radius=0.1]
		          (0,2) circle [radius=0.1];
		
		\begin{scope}[shift={(5,0)}]
    		\draw[postaction={decorate}] (0,0) -- node[below] {$2 \cdot 1$} (2,0);
    		\draw[postaction={decorate}] (0,0) -- node[left] {$2 \cdot 2$} (0,2);
    		\draw[postaction={decorate}] (2,2) -- node[right] {$3 \cdot 1$} (2,0);
    		\draw[postaction={decorate}] (2,2) -- node[above] {$3 \cdot 2$} (0,2);
    		
    		\draw[fill=white] (0,0) circle [radius=0.1] node[above right] {2}
    		                  (2,2) circle [radius=0.1] node[below left] {3};
    		\filldraw (2,0) circle [radius=0.1] node[above left] {1}
    		          (0,2) circle [radius=0.1] node[below right] {2};
		\end{scope}
	\end{tikzpicture}
    \caption{A class whose edge labels are induced by vertex labels. The bipartition of the vertices is shown in black and white.}
    \label{figEx}
\end{figure}

Recall that a class $\alpha \in H_2(A_\G)$ is visually represented by integer labels on the edges of the graph~$\G$ equipped with an orientation. 
The class $\alpha = \overline v \otimes \overline w$ has the special property that there is an integer labelling of the vertices in $\G$ that induces the edge labelling in the following way: the edge label is given by multiplying the vertex labels on the edge's incident vertices; the vertex labels are precisely the coordinates of~$\overline v \in H_1(F_n) \cong \Z^n$ and~$\overline w \in H_1(F_m) \cong \Z^m$ with respect to the canonical bases. See Figure~\ref{figEx} for an illustration. 

Surprisingly, this simple example allows us to compute the minimal genus of any class.

\begin{thrm}\label{thrm - inequality is equality complete bipartite}
Let $\G$ be a complete bipartite graph and $\al \in H_2(A_\G) \cong H_1(F_n) \otimes H_1(F_m)$. Then the minimal genus $\gen(\al)$ is equal to the cap bound $ \frac{1}{2} \rk(M_\al)$ and the minimal number of pure tensors, i.e.~elements of the form~$\overline v \otimes \overline w$, needed to represent~$\al$. Furthermore, the minimal genus is always realised by a disjoint union of tori.

\begin{proof} Let~$\G$ be a complete bipartite graph with vertices $\{v_1, \ldots, v_n \} \sqcup \{ w_1, \ldots, w_m \}$ and an orientation given by $E^{or}(\G)=\{v_1, \ldots, v_n \} \times \{ w_1, \ldots, w_m \}$. Suppose $\alpha \in H_2(A_\G)$ is an arbitrary class. We want to show that $2 \, \gen(\alpha) = \rk(M_\alpha)$. This holds automatically if~$\al$ is trivial, so we may assume it is nontrivial. Recall that the connection matrix is
\[
M_\al = 
\begin{pNiceArray}{ccc|ccc}[first-row, first-col]
 & v_1 & \Cdots & v_n & w_1 & \Cdots & w_m \\
v_1 & \Block{3-3}<\Large>{0} &  &  & l(v_1,w_1) & \cdots & l(v_1,w_m) \\
\Vdots & &  & & & \ddots & \\
v_n &  &  &  & l(v_n,w_1) & \cdots & l(v_n,w_m) \\
\hline
w_1 & -l(v_1,w_1) & \cdots & -l(v_n,w_1) & \Block{3-3}<\Large>{0} &  & \\
\Vdots & & \ddots & & &  & \\
w_m & -l(v_1,w_m) & \cdots & -l(v_n,w_m) &  &  & 
\end{pNiceArray}.
\]
Note that $\rk(M_\alpha)$ is twice the rank of the submatrix given by the first $n$ rows and the last $m$ columns. We start with the case $\rk(M_\alpha) = 2$. In this case, the rows of the submatrix span a cyclic subgroup in $\Z^m$. This implies there is a nontrivial integral row vector $(z_j)_{j=1}^m$ and integral multipliers $(c_i)_{i=1}^n$ such that $c_i \cdot z_j = l(v_i, w_j)$ for $1 \le i \le n$ and $1 \le j \le m$. In other words, $\alpha = \overline x \otimes \overline y$ if we set~$\overline x = \sum_{i=1}^n c_i \overline v_i$ and~$\overline y = \sum_{j=1}^m z_j \overline w_j$ in~$H_1(F_n)$ and $H_1(F_m)$ respectively. By Lemma~\ref{Thorben}, this implies~$\alpha$ is representable by a torus and $2 \, \gen(\alpha) = 2 = \rk(M_\alpha)$.

More generally, suppose $\rk(M_\alpha) = 2k$. Then there are $k$ linearly independent integral row $m$-vectors $\overline z_1 = (z_{1j})_{j=1}^m, \ldots,\,\overline z_k = (z_{kj})_{j=1}^m$ and $k$ integral $n$-vectors $\overline c_1 = (c_{1i})_{i=1}^n, \ldots,\,\overline c_k = (c_{ki})_{i=1}^n$ such that $ c_{1i} \cdot z_{1j} + \cdots +  c_{ki} \cdot z_{kj} = l(v_i, w_j)$ for $1 \le i \le n$ and $1 \le j \le m$. As before, this means $\alpha = \sum_{l=1}^k \overline x_l \otimes \overline y_l$ where we set~$\overline x_l = \sum_{i=1}^n c_{li} \overline v_i$ and~$\overline y_l = \sum_{j=1}^m z_{lj} \overline w_j$ for~$l = 1, \ldots, k$. But each summand has a torus representative by Lemma~\ref{Thorben}, and so $2 \, \gen(\alpha) = 2k = \rk(M_\alpha)$.
\end{proof}
\end{thrm}

In particular, this theorem characterises classes representable by a torus as those considered in Lemma~\ref{Thorben}:

\begin{cor}\label{cor - torusrep} Let~$\G$ be a complete bipartite graph and $A_\G \cong F_n \times F_m$. A class $\alpha \in H_2(A_\G)$ is representable by a torus if and only if $\alpha = \overline v \otimes \overline w$ for some $\overline v \in H_1(F_n)$ and $\overline w \in H_1(F_m)$.
\begin{proof}
The reverse direction is precisely the construction in Lemma~\ref{Thorben}. For the forward direction, suppose $\alpha \in H_2(A_\G)$ is nontrivial and $\gen(\alpha) = 1$. Then, by Theorem~\ref{thrm - inequality is equality complete bipartite}, this implies $\alpha = \overline v \otimes \overline w$ for some $\overline v \in H_1(F_n)$ and $\overline w \in H_1(F_m)$.
\end{proof}
\end{cor}

\begin{defn}
We call a simple graph a \emph{star} if it consists of one vertex of valence~$n$,~$n$ edges, and~$n$ leaves. In other words, it is a complete bipartite graph with partition of the vertices $V(\G) = \{ v_1, \ldots, v_n \} \sqcup \{ w \}$.
\end{defn}

The following lemma is another immediate consequence of Lemma~\ref{Thorben}.

\begin{lem}\label{lemma-star}
Let~$\G$ be any simple graph and~$\alpha \in H_2(A_{\G})$. If~$\supp(\alpha)$ is a star, then~$\alpha$ is representable by a torus.
\begin{proof}
Fix an orientation of~$\G$. If~$\G'=\supp(\alpha)$ is a star with vertex set $\{ v_1, \ldots, v_n \} \sqcup \{ w \}$, then we can induce the appropriate edge labels of~$\alpha$ from vertex labels: label the vertex~$w$ with~$1$ and the vertex~$v_i$ with the same label as the edge~$\{w,v_i\}$. 
Let~$\alpha' \in H_2(A_{\G'})$ be the class whose edge labels match those of $\alpha$. By the discussion after Lemma~\ref{Thorben}, the class~$\alpha'$ is representable by a torus. The inclusion $\G' \subseteq \G$ induces $\iota_*:H_2(A_{\G'}) \to H_2(A_{\G})$ with $\iota_*(\alpha') = \alpha$. Therefore, $\alpha$ is representable by a torus.
\end{proof}
\end{lem}

\subsection{Trees}\label{subsect - tree}
In this subsection, we compute the minimal genus in the case where~$\G$ is a tree. Recall that a tree is a simple graph with no cycles. Some of our results hold for general graphs~$\G$ and so we will always make explicit our assumptions on the graph.

The cap bound allows us to constrain the options for~$\supp(\alpha)$ when~$\alpha$ is representable by a torus.

\begin{lem}\label{lem - two edges form square in sup(alpha)}
Let $\G$ be any simple graph and~$\alpha\in H_2(A_\G)$. If~$\rk(M_\alpha)=2$, then any two edges with distinct vertices~$\{v_1,v_2\}$ and~$\{w_1,w_2\}$ in~$\supp(\alpha)$ form two sides of a square in~$\supp(\al)$.
\begin{proof}
Suppose~$\alpha \in H_2(A_\G)$ satisfies~$\rk(M_\alpha)=2$, and suppose $\{v_1,v_2\}$ and~$\{w_1,w_2\}$ are edges in~$\supp(\alpha)$ with~$l(v_1,v_2)=\lambda\neq 0$ and~$l(w_1,w_2)=\beta\neq 0$. For brevity, set~$l^i_j=l(v_i,w_j)$. Then the matrix~$M_\alpha$ has the following submatrix:

\[
M_\al |_{\langle v_1,v_2,w_1,w_2 \rangle}=
\kbordermatrix{
 & v_1 & v_2 & w_1 & w_2 \\
v_1 & 0 & \lambda & l^1_1 & l^1_2 \\
v_2 & -\lambda & 0 & l^2_1 & l^2_2  \\
w_1 & -l^1_1 & -l^2_1& 0& \beta  \\
w_ 2 & -l^1_2 & -l^2_2 &  -\beta&0}.
\]

If the two edges are disjoint, all~$l^i_j=0$ and the submatrix has rank 4. So~$\rk(M_\al)\geq 4$ and this contradicts $\rk(M_\alpha)=2$.
The other option if~$\{v_1,v_2\}$ and $\{w_1,w_2\}$ do not make a square is that~$l^1_1=l^1_2=0$ but one or both of~$l^2_1$ and~$l^2_2$ are non-zero. The matrix becomes

\[
M_\al |_{\langle v_1,v_2,w_1,w_2 \rangle}=
\kbordermatrix{
 & v_1 & v_2 & w_1 & w_2 \\
v_1 & 0 & \lambda & 0 & 0 \\
v_2 & -\lambda & 0 & l^2_1 & l^2_2  \\
w_1 & 0 & -l^2_1& 0& \beta  \\
w_ 2 & 0 & -l^2_2 &  -\beta&0},
\]
which again has rank $4$ and contradicts the hypothesis.
\end{proof}
\end{lem}

As a corollary, any class representable by a torus has a connected support.

\begin{cor}\label{cor - supp connected for torus rep}
Let $\G$ be a simple graph and~$\alpha\in H_2(A_\G)$. If~$\rk(M_\al) = 2$, then the support~$\supp(\alpha)$ is connected.
\begin{proof}
We prove the contrapositive, assuming~$\al \in H_2(A_\G)$ is nontrivial. Suppose the support is disconnected and consider two edges taken from two connected components. Then these edges do not form a square. So~$\rk(M_\al) \ge 4$ by Lemma~\ref{lem - two edges form square in sup(alpha)}.
\end{proof}
\comment{
\begin{proof}
Suppose that~$\supp(\alpha)$ is disconnected. Since~$\supp(\alpha)$ is a union of labelled edges,it follows that there are two edges~$(v_1,v_2)$ and~$(w_1,w_2)$ in~$\supp(\alpha)$ such that~$v_i$ is not connected to~$w_j$ for~$i,j \in \{1,2\}$. Let~$l(v_1,v_2)=\lambda$ and~$l(w_1,w_2)=\beta$. The matrix~$M_\al$ has the following submatrix.
\[
M_\al |_{<v_1,v_2,w_1,w_2>}=
kbordermatrix{
 & v_1 & v_2 & w_1 & w_2 \\
v_1 & 0 & \lambda & 0 & 0 \\
v_2 & -\lambda & 0 & 0 & 0 \\
w_1 & 0 & 0& 0& \beta  \\
w_ 2 & 0 & 0 &  -\beta&0}
\]
Since this submatrix has rank 4, it follows~$\rk(M_\al) \geq 4$, which bounds the genus below by 2 using Equation~\ref{eq - cap bound}. It follows~$\al$ is not representable by a torus. 
\end{proof}
}
\end{cor}

\begin{prop}\label{prop- star in tree is torus}
Let~$\G$ be a tree. A class~$\alpha\in H_2(A_\G)$ is representable by a torus if and only if $\supp(\alpha)$ is a complete bipartite graph (or, alternatively, a star).
\begin{proof}

The statement is vacuously true when~$\al$ is trivial. So we may assume that~$\al$ is nontrivial. If $\supp(\alpha)$ is a complete bipartite graph, then~$\supp(\alpha)$ is a star since~$\G$ is a tree. From Lemma~\ref{lemma-star}, it follows that~$\alpha$ is representable by a torus.

Suppose, conversely, that $\al$ is representable by a torus but $\supp(\alpha)$ is not a complete bipartite graph, i.e.~not a star. Then since~$\al$ is representable by a torus,~$\supp(\al)$ is connected by the cap product inequality (Equation~\ref{eq - cap bound}) and Corollary~\ref{cor - supp connected for torus rep}. Moreover, as~$\G$ is a tree, there is a subgraph of~$\supp(\al)$ of the following form:

\begin{center}
    	\begin{tikzpicture}[thick, scale=.7]
			
			\begin{scope}[scale=1.2,shift={(-5,3)},decoration={
        markings,
        mark=at position 0.5 with {\arrow{>}}}]
			
			\draw[rotate =0] (0,0)--(1,0);
			\draw[rotate =0] (1,0)--(2,0);
			\draw[rotate =60] (0,0)--(1,0);
			\draw[rotate =120] (0,0)--(1,0);
			\draw[rotate =220] (0,0)--(1,0);
			\draw[rotate =180] (0,0)--(1,0);
			
			\draw[rotate =310] (.8,0)  circle [radius=0.01];
			\draw[rotate =300] (.805,0)  circle [radius=0.01];
			\draw[rotate =290] (.8,0)  circle [radius=0.01];

			\filldraw (0,0) circle [radius=0.1] node[below] {$w$};
			\filldraw (1,0) circle [radius=0.1] node[above] {$v_k$};
			\filldraw (2,0) circle [radius=0.1] node[above] {$z$};
			\filldraw (-1,0) circle [radius=0.1] node[left] {$v_3$};
			\filldraw[rotate =60] (1,0) circle [radius=0.1] node[above] {$v_1$};
			\filldraw[rotate =120] (1,0) circle [radius=0.1] node[above] {$v_2$};
			\filldraw[rotate =220]  (1,0) circle [radius=0.1] node[below] {$v_4$};
			\end{scope}
		\end{tikzpicture}
\end{center}

\noindent where~$k$ is at least~$2$.
 The edges~$\{w, v_1\}$ and~$\{v_k,z\}$ cannot form two sides of a square in $\supp(\al)$ since~$\G$ is a tree. By Lemma~\ref{lem - two edges form square in sup(alpha)}, we have~$\rk(M_\al) \ge 4$. This contradicts the cap product inequality (Equation~\ref{eq - cap bound}) and the assumption~$\al$ was representable by a torus; therefore, the support~$\supp(\al)$ is a complete bipartite graph.
\end{proof}
\end{prop}

\begin{defn}
A \emph{star covering} of an integer labelled oriented simple graph~$\G$ with labels~$l(v,w)$ is a finite collection of integer labelled oriented star graphs $\{S_1,\ldots, S_k\}$ such that:
\begin{enumerate}
    \item for all~$i$, $S_i$ is a subgraph of $\G$, and the orientation on~$S_i$ is induced by the orientation on~$\G$; and
    \item let~$l_i(v,w)$ be the label of~$(v,w)$ if it lies in~$S_i$, and zero otherwise; then for all oriented edges~$(v,w)\in E^{or}(\G)$, we require $\sum_{i=1}^k l_i(v,w)=l(v,w)$. 
\end{enumerate}
We denote by~$\stc(\G)$ the minimal~$k$ for which there exists such a star covering of~$\G$. We also write $\stc(\al)$ for $\stc(\supp(\al))$, where we now consider~$\supp(\al)$ as a labelled oriented graph. Note that~$\stc(\al)$ is independent of the orientation chosen to depict~$\al$.
\end{defn}

\begin{remark}
In the literature an unlabelled version of this cover is often called a \emph{vertex cover}. Since our emphasis lies on homology classes described by labelled oriented edges we use the term star cover instead.
\end{remark}

\begin{lem} \label{lem - minimal star covering disjoint}
Given a star covering of a labelled oriented graph~$\G$, there is a covering~$\{S_1,\ldots S_k\}$ with the same cardinality so that for all~$i\neq j$,~$S_i$ and~$S_j$ have disjoint edge sets (their vertex sets may intersect non-trivially). 
\begin{proof}
Suppose multiple stars~$\{S_i\}_{i\in I}$ have a common oriented edge~$(v,w)$ for some $I\subseteq \{1,\ldots k\}$ and that~$S_i$ has label~$l_i(v,w)\neq 0$ for~$i\in I$. Choose one~$j \in I$ and remove the edge~$(v,w)$ from every star $S_i$ with $j \neq i \in I$. Following this, change the label on~$S_j$ to~$l_j(v,w) = l(v,w)$. Repeat this process for all edges in~$\G$ which are common edges between stars.
\end{proof}
\end{lem}

\begin{prop}\label{prop - star covering bound}
Let~$\G$ be any oriented simple graph and~$\al \in H_2(A_\G)$. The minimal cardinality of a star covering of~$\supp(\alpha)$ bounds the minimal genus from above, i.e.
\[\stc(\alpha)\geq \gen(\alpha).\]
\end{prop}
Note that, a priori, an upper bound on the minimal genus for a second homology class of a RAAG is given by the number of edges in the support. This proposition provides a substantial improvement to this bound.
\begin{proof}
Given a star covering $\{S_1,\ldots S_k\}$ of minimal cardinality of~$\supp(\alpha)$ (so~$k=\stc(\al)$), by Lemma~\ref{lem - minimal star covering disjoint}, we can assume the edge sets of the stars are disjoint, i.e.~we can label the edges of the stars with the edge labels of~$\alpha$. For~$1\leq i \leq k$, we associate to the star~$S_i$ a class~$s_i \in H_2(A_\G)$ with labelled support~$S_i$ such that
\[
\alpha=\sum_{i=1}^k s_i.
\]
By Lemma~\ref{lemma-star}, we know that~$g(s_i)=1$ for~$1\leq i\leq k$ and so it follows from subadditivity (Lemma~\ref{lem-subadditive}) that
\[\gen(\alpha)\leq k = \stc(\alpha).\qedhere
\]
\end{proof}

\begin{prop}\label{prop - stc equals cap}
Let $\G$ be an oriented tree and $\al \in H_2(A_\G)$.
Then the cap bound $\frac{1}{2} \rk(M_\al)$ is equal to~$\stc(\al)$.
\begin{proof}
From Proposition~\ref{prop - star covering bound}, we know that $\stc(\alpha)\geq \gen(\alpha)$, and by Equation \ref{eq - cap bound}, this implies that $\stc(\alpha)\geq \frac{1}{2} \rk(M_\al)$. We now show the reverse inequality. The result is obvious if~$\al$ is trivial, so we may assume it is nontrivial.

First of all, we  assume that $\supp(\al)$ is connected. Let~$\{S_1,\ldots S_k\}$ be a star covering of minimal cardinality on~$\supp(\al)$ such that all stars have disjoint edge sets (Lemma~\ref{lem - minimal star covering disjoint}). We proceed by induction on~$k = \stc(\al)$. If~$k=1$, then~$\supp(\alpha)$ is a star and, by Proposition~\ref{prop- star in tree is torus},~$\gen(\alpha)=1$. This implies $\frac{1}{2} \rk(M_\al)=1=\stc(\alpha)$. Assume that~$k>1$ and that the theorem holds for classes~$\al'\in H_2(A_\G)$ with~$\stc(\alpha')\leq k-1$. Then, since~$\G$ is a tree,~$\supp(\alpha)$ has the following form
\begin{center}
    	\begin{tikzpicture}[thick, scale=.7]
	    	\begin{scope}[scale=1.2,shift={(-5,3)},decoration={markings, mark=at position 0.5 with {\arrow{>}}}]
        
			\draw[rotate =0, postaction=decorate] (0,0)--(1,0);
			\draw[rotate =60, postaction=decorate] (0,0)--(1,0);
			\draw[rotate =120, postaction=decorate] (0,0)--(1,0);
			\draw[rotate =220, postaction=decorate] (0,0)--(1,0);
			\draw[rotate =180, postaction=decorate] (0,0)--(1,0);

			\draw (1,0)--(1.3,0.25) (1,0)--(1.3,-.25) (1,0)--(1.4,0);
			
			\draw[rotate =310] (.8,0)  circle [radius=0.01];
			\draw[rotate =300] (.805,0)  circle [radius=0.01];
			\draw[rotate =290] (.8,0)  circle [radius=0.01];
			
			\draw (1,-.85)--(1,.85)--(5,.85)--(5,-.85)--(1,-.85);
			
			\filldraw (0,0) circle [radius=0.1] node[below] {$w$};
			\filldraw (1,0) circle [radius=0.1] node[above, xshift=-1.1ex] {$v_p$};
			\filldraw (-1,0) circle [radius=0.1] node[left] {$v_3$};
			\filldraw[rotate =60] (1,0) circle [radius=0.1] node[above] {$v_1$};
			\filldraw[rotate =120] (1,0) circle [radius=0.1] node[above] {$v_2$};
			\filldraw[rotate =220]  (1,0) circle [radius=0.1] node[below] {$v_4$};
			
			\draw (3,0) node {$\G'$};
			\end{scope}
		\end{tikzpicture}
\end{center}
where the star pictured is~$S_k$,~$v_1$ is a leaf of~$\supp(\alpha)$, and the stars~$\{S_1,\ldots, S_{k-1}\}$ form a star covering of minimal cardinality for~$\G'=\supp(\al)\backslash S_k$. Then there is a class~$\al' \in H_2(A_\G)$ with~$\supp(\al')=\G'$ and~$\stc(\al')=k-1$. By the induction hypothesis,~$\frac{1}{2} \rk(M_{\al'})=k-1$.

Let~$u$ be a leaf of~$\supp(\al)$ and~$S_k$. If~$S_k$ only has one edge~$\{u,v\}$, then there is another star~$S_j$ in the star cover containing an edge~$\{v,z\}$ for some~$z\neq v$ (because $\supp(\al)$ is connected and $k>1$). Remove the edge~$\{v,z\}$ from~$S_j$ and add it to~$S_k$ along with its orientation and label. This now gives a star covering of the same cardinality such that~$S_k$ is as pictured above with~$p\geq 2$. If~$S_k$ has two or more edges, we automatically have~$p \geq 2$ and no modifications are needed. Either way, we assume~$S_k$ was chosen such that~$p\geq 2$. 

Consider the connection matrix~$M_\alpha$ restricted to~$v_1, w$, and the vertices in~$\G'$, and let~$l(v_1,w)=\lambda$ and~$l(w,z_i)=\beta_i$ where~$\{z_1,\ldots ,z_r\}$ is the vertex set of~$\G'$.
\[
M_\al|_{\langle v_1,w, V(\G') \rangle}=
\begin{pNiceArray}{cc|ccc}[first-row, first-col]
&v_1 & w & \Block{1-3}{V(\G')}&&\\
v_1 & 0 & \lambda & 0 &\Cdots & 0\\
w& -\lambda &0 & \beta_1 & \Cdots & \beta_r\\
\hline 
 \Block{3-1}{\rotate V(\G')}&0& -\beta_1 &\Block{3-3}<\Large>{M_{\al'}}&&\\
& \Vdots & \Vdots & &&\\
&0 &-\beta_r & &&
\end{pNiceArray}
\]
Observe that~$\rk(M_\al)\geq 2 + \rk(M_{\al'})$. Since~$\frac{1}{2} \rk(M_{\al'})=k-1$, it follows that $\frac{1}{2} \rk(M_\al)\geq k=\stc(\al)$ and this shows the reverse inequality we required. So~$\frac{1}{2} \rk(M_\al)=\stc(\al)$.

Now suppose that~$\supp(\al)$ is disconnected. Then \[
\alpha=\sum_{j=1}^l \beta_j
\]
for some~$\beta_j\in H_2(A_\G)$ such that~$\supp(\beta_j)$ is connected, and~$l$ is the number of connected components of~$\supp(\alpha)$. Then since the~$\beta_j$ all have disjoint support, $M_\al$ is a block diagonal matrix with blocks~$M_{\beta_j}$ for~$1\leq j \leq l$. It follows from computing the rank that
\[
\frac{1}{2}\rk(M_\al) =\frac{1}{2} \sum_{j=1}^l \rk(M_{\beta_j}) = \sum_{j=1}^l \stc(\beta_j)
\]
where the final equality comes from the first part of this proof applied to each~$\beta_j$. Finally, we note that since each~$\beta_j$ corresponds to a connected component of~$\supp(\alpha)$, any star covering of~$\alpha$ is a union of star coverings of the~$\beta_j$ and vice versa. So~$\frac{1}{2} \rk(M_\al) = \sum_{j=1}^l \stc(\beta_j)=\stc(\alpha)$ as required.
\end{proof}
\end{prop}

\begin{thrm}\label{thrm - inequality is equality trees}
Let $\G$ be a tree and $\al \in H_2(A_\G)$. Then the minimal genus $\gen(\al)$ is equal to the cap bound $\frac{1}{2} \rk(M_\al)$ and~$\stc(\al)$. Furthermore, the minimal genus is always realised by a disjoint union of tori.

\begin{proof}
Fix an orientation on the tree~$\G$. From Proposition~\ref{prop - star covering bound}, we know that $\stc(\alpha)\geq \gen(\alpha)$. On the other hand, from a combination of Proposition~\ref{prop - stc equals cap} and Equation \ref{eq - cap bound}, we have
\[\stc(\al)=\frac{1}{2} \rk(M_\al)\leq \gen(\al).\]
Putting these two inequalities together gives $\stc(\al)=\gen(\al)$. Recall that each star in the star covering of minimal cardinality is representable by a torus by Lemma~\ref{lemma-star} and the disjoint union of these tori has genus $\stc(\al)=\gen(\al)$. Thus the minimal genus is always realised by a disjoint union of tori.
\end{proof}
\end{thrm}

\section{Classes representable by a torus}\label{section - one torus}
In this section, we investigate the relationship between the support of a class and its minimal genus. We restrict ourselves to the case where the minimal genus is one, i.e.~the class is representable by a torus.

Recall that a \emph{complete~$n$-partite graph}~$\G$ has a partition of the vertices $V(\G)=\sqcup_{i=1}^n X_i$ such that for every $1\leq i \leq n$:
\begin{itemize}
    \item all $x_i\in X_i$ are joined by an edge to every vertex in~$V(\G)\backslash X_i$
    \item there exist no edges between any pair of vertices in~$X_i$.
\end{itemize}
Given an arbitrary graph~$\G$, we also use the notion of a \emph{maximal complete $m$-partite subgraph}. This is a full subgraph~$Y$ of~$\G$ such that~$Y$ is complete $m$-partite for some~$m$, and~$V(Y)$ is maximal over all full, complete~$k$-partite subgraphs of~$\G$ for any~$k$. Note that such a maximal complete $m$-partite subgraph is not unique. Recall that a full subgraph is one which inherits all edges that are present in~$\G$, i.e.~if~$v, w \in V(Y)$ and~$\{v,w\} \in E(\G)$, then~$\{v,w\}\in E(Y)$.

\begin{prop}
\label{prop - rank 2 implies n partite}
Let~$\G$ be a simple graph and~$\alpha \in H_2(A_\G)$. If~$\rk(M_\alpha) = 2$, then~$\supp(\alpha)$ is a complete $n$-partite graph for some positive integer $n$.

\begin{proof}
Since~$\rk(M_\alpha)=2$, the support~$\supp(\al)$ is nonempty and it follows from Corollary~\ref{cor - supp connected for torus rep} that~$\supp(\al)$ is connected.

We assume that~$\supp(\al)$ is not complete~$n$-partite for any integer~$n \ge 1$ and work towards a contradiction. Consider a maximal complete $m$-partite subgraph~$Y$ of~$\supp(\al)$ with~$V(Y)=\sqcup_{i=1}^m X_i$. Then~$Y\neq\supp(\alpha)$ by the assumption and, since~$\supp(\al)$ is connected, there exists a vertex~$w$ in~$\supp(\al)\backslash Y$ connected by an edge of~$\supp(\al)$ to some vertex in~$Y$. 

{\bf Claim:} If~$w$ is connected by an edge of~$\supp(\al)$ to some vertex in~$X_j$, then~$w$ is connected by an edge of~$\supp(\al)$ to every vertex in~$X_j$.

\emph{Proof of claim.} Suppose that~$w$ is connected by an edge of~$\supp(\al)$ to~$x \in X_j$ but not to~$y\in X_j$. We consider the connection matrix for~$\alpha$ restricted to~$w,x,y,$ and~$V(Y)\backslash X_j$. Let~$l(w,x)=\lambda \neq 0$ and denote the labels between~$y$ and~$V(Y)\backslash X_j$ by~$\mu_i \neq 0$ for~$1\leq i \leq r$. The labels~$\mu_i$ are non-zero since~$Y$ is a complete~$m$-partite subgraph of~$\supp(\al)$. (We mark unimportant entries with stars.)
\[
M_\al|_{\langle w,x,y, V(Y)\backslash X_j \rangle}=
\begin{pNiceArray}{ccc|ccc}[first-row, first-col]
&w & x & y & \Block{1-3}{V(Y)\backslash X_j }&&\\
w & 0 & \lambda &0& * & \Cdots& *\\
x& -\lambda &0 & 0 &* & \Cdots & * \\
y& 0 &0 & 0& \mu_1 & \Cdots & \mu_r\\
\hline 
 \Block{3-1}{\rotate V(Y)\backslash X_j }& *& *&-\mu_1 &\Block{3-3}<\Large>{M_{\al}|_{\langle V(Y)\backslash X_j \rangle}}&&\\
& \Vdots &\Vdots & \Vdots &&&\\
& *&*&-\mu_r & &&
\end{pNiceArray}
\]

Since this is a submatrix of~$M_\al$, and the~$\mu_i$ are all non-zero, it follows that~$\rk(M_\al)\geq 3$, which contradicts the hypothesis that~$\rk(M_\al) = 2$. This proves the claim.

Thus for each~$X_j\subset V(Y)$,~$w$ is connected (in~$\supp(\al)$) either to all or none of the vertices in~$X_j$. For the contradiction, we now rule out all possible ways~$w$ may be connected to~$Y$ by an edge of~$\supp(\al)$. 

First, we note that $w$ is not connected by an edge of~$\supp(\al)$ to all~$v\in Y$ because then the full subgraph of~$\supp(\al)$ spanned by~$w$ and $V(Y)$ would be a complete $(m+1)$-partite subgraph of~$\supp(\al)$ with one more vertex than~$Y$, and this contradicts the maximality of~$Y$. 

For the same reason,~$w$ cannot be connected by an edge of~$\supp(\al)$ to all vertices in~$V\backslash X_j$ for some~$j$, and disconnected from~$X_j$: replacing~$X_j$ with $\{ w \} \cup X_j$ would result in the full subgraph of~$\supp{\al}$ spanned by~$w$ and $V(Y)$ being a complete $m$-partite subgraph of~$\supp(\al)$ with one more vertex than~$Y$.

The final case to check is when~$w$ is not connected to~$l$ of the~$X_j$, but is connected to~$(m-l)$ of the~$X_j$, for~$1<l<m$. Without loss of generality, suppose~$w$  is disconnected from all vertices in~$\sqcup_{j=1}^l X_j$ and connected by an edge to all vertices in~$\sqcup_{j=l+1}^m X_j$. Then the submatrix~$M_\al|_{\langle w, V(Y) \rangle}$ is given by:
\[
M_\al|_{\langle w, V(Y) \rangle}=
\begin{pNiceArray}{c|ccc|ccc}[first-row, first-col]
&w & X_1 & \Cdots & X_l& X_{l+1}&\Cdots&X_m\\
w & 0 & 0&\cdots& 0 & \mu_{l+1} &\cdots& \mu_m \\
\hline
X_1& 0 &\Block{3-3}<\Large>{M_{\al}|_{\langle X_{1},\ldots , X_l \rangle}} &  & & \Block{3-3}<\Large>{N} &  \\
\Vdots& \vdots  & & &&&  & \\
X_l& 0 & & & & &  & \\
\hline 
 X_{l+1}& -\mu_{l+1}&  \Block{3-3}<\Large>{-N}& &&\Block{3-3}<\Large>{M_{\al}|_{\langle X_{l+1},\ldots , X_m \rangle}}& &\\
\vdots & \Vdots & &  &&&\\
X_m& -\mu_m&& & &&
\end{pNiceArray}
\]
where the matrix~$N$ has no zero entries and~$\mu_j \neq 0$ for~$l+1\leq j \leq m$. Then, since~$l$ is greater than~$1$, $\rk(M_{\al}|_{\langle X_{1},\ldots , X_l \rangle})\geq 2$ and thus~$\rk(M_\al|_{\langle w, V(Y) \rangle})\geq 3$ (since the~$\mu_j$ are non-zero), which contradicts~$\rk(M_\al)=2$.
\end{proof}

\end{prop}

\begin{prop}\label{prop - n-partitite + rank 2 implies torus}
Let~$\G$ be a complete n-partite graph and~$\alpha \in H_2(A_\G)$. If $\rk(M_\alpha) = 2$, then~$\alpha$ is representable by a torus.

\begin{proof}
The abelianisation of~$A_\G$ is~$\Z^{V(\G)}$, and the map from~$A_\G$ to its abelianisation corresponds to embedding~$\G$ into the complete graph on the vertex set~$V(\G)$. Let $\alpha^{ab} \in H_2(\Z^{V(\G)})$ denote the image of~$\al$ under the homomorphism on second homology induced by the abelianisation. Then the connection matrices are equal, i.e.~$M_{\alpha^{ab}} = M_\alpha$. 

Recall that~$H_2(\Z^{V(\G)})$ is identified with $\bigwedge^2 \Z^{V(\G)}$. Since~$\rk(M_\alpha) = 2$, it follows from Proposition~\ref{prop - minimal number of elementary wedges} that the class~$\alpha^{ab}$ is an elementary wedge $a \wedge b$ for some vectors~$a,b \in \Z^{V(\G)}$.

Using the partition $V(\G)=\sqcup_{i=1}^n X_i$, write $a=\sum_{i=1}^n a_i,~b=\sum_{i=1}^n b_i$ where each~$a_i, b_i$ lie in the subgroup generated by $X_i$.
The partition of~$V(\G)$ induces a block decomposition of~$M_\alpha$ with diagonal blocks~$M_\alpha|_{\langle X_i\rangle}$; these diagonal blocks are necessarily zero matrices since there are no edges in~$\G$ between vertices in~$X_i$. Since~$M_\al=M_{\al^{ab}}$, the blocks~$M_\alpha|_{\langle X_i\rangle}$ correspond to the elementary wedges~$a_i \wedge b_i$ under the identification $\alpha^{ab}=a \wedge b$, and it follows that~$a_i \wedge b_i$ is trivial.
This implies~$a_i$ and~$b_i$ are linearly dependent, i.e.~there are vectors~$v_i$ and multipliers~$r_i, s_i \in \Z$ such that~$a_i = r_i v_i,~b_i = s_i v_i$ for~$i= 1, \ldots, n$.

Choose lifts~$\nu_i \in A_\G$ of the vectors~$v_i \in \Z^{V(\G)}$. As each element $\nu_i$ lies in the subgroup generated by~$X_i$, the set of elements~$\nu_i$ pairwise commute and thus determine a homomorphism $\tau: \Z^n \to A_\G$ that maps the standard basis $e_i$ to $\nu_i$. Set $c = \sum_{i=1}^n r_i e_i$ and $d = \sum_{i=1}^n s_i e_i$; using the Pontryagin ring structures on homology groups of $\Z^n$ and $\Z^{V(\G)}$, we see that $\tau_*(c \wedge d)^{ab} = a \wedge b$, where $\tau_*\colon H_2(\Z^n) \to H_2(A_\G)$ is the homomorphism induced by $\tau$. So $\alpha = \tau_*\left(c \wedge d\right)$ since the abelianisation map is injective on second homology. By Proposition~\ref{prop - n torus cap equality}, the elementary wedge~$c \wedge d$, and hence its image $\alpha$, is representable by a torus.
\end{proof}

\end{prop}

Putting the results of this section together gives:
\setcounter{abcthm}{3}
\begin{abcthm}\label{abcthm - n partite torus}
Let~$\G$ be any simple graph. Then a nontrivial second homology class $\alpha \in H_2(A_\G)$ is representable by a torus if and only if $\rk(M_\alpha) = 2$. Furthermore, the support of such a homology class is a complete n-partite graph.
\begin{proof}
If~$\al \in H_2(A_\G)$ is nontrivial, then automatically~$\rk(M_\al) \ge 2$. If we also assume~$\al$ is representable by a torus, then $\rk(M_\al) = 2$ by the cap product inequality (Equation~\ref{eq - cap bound}). Conversely, if~$\rk(M_\al) = 2$, then~$\al$ is nontrivial and, by Proposition~\ref{prop - rank 2 implies n partite}, its support~$\supp(\al)$ is a complete $n$-partite graph for some positive integer~$n$. Restrict the graph~$\G$ if necessary and assume~$\supp(\al) = \G$. By Proposition~\ref{prop - n-partitite + rank 2 implies torus}, the class~$\al$ is representable by a torus.
\end{proof}
\end{abcthm}

\section{An example answering Question~\ref{question - realised by dj tori}}\label{section - question 2}
In this final section, we answer the following question:
\begin{question}  Does every class in the second homology of a RAAG have a minimal genus representative that is a disjoint union of tori?
\end{question}
To describe a second homology class representative in general, we use \emph{Van Kampen diagrams} (see~\cite{VKDiagrams} for a more thorough introduction). 
These diagrams encode a cellular structure of a compact surface~$\Sigma$ and information on how the cells are mapped to a CW space~$X$: overall this gives a cellular map~$f\colon \Sigma \to X$. In our case, we want a map from a surface $\Sigma$ to $\Sal_\G$, so the Van Kampen diagram is a tesselation of $\Sigma$ by squares with oriented edges labelled by vertices of $\G$ such that opposite edges of the squares have the same labelling and orientation  and for each square either:
\begin{enumerate}
    \item all edges are labelled by the same generator, or
    \item edges are labelled by two generators $v$ and $w$ such that $\{v,w\}\in E(\G)$
\end{enumerate}
On the $1$-skeleton, the corresponding map $f\colon \Sigma \to \Sal_\G$ is given by mapping each vertex to the single vertex $x_0$ of $\Sal_\G$ and every edge to the 1-cube (or circle) in the Salvetti complex corresponding to its label. Since the boundary of a square is mapped to the commutator of its labels, the map on the $1$-skeleton can be extended to the 2-skeleton if and only if said commutators vanish. In our case this holds since either the commutator is trivially trivial (case~1 above) or there is an edge in~$\G$ so the commutator is trivial in~$A_\G$ (case~2 above). We note that such an extension is unique up to homotopy as $\pi_2(\Sal_\G)$ is trivial. 
Hence a Van Kampen diagram encodes a map $f\colon \Sigma\to \Sal_\G$ and it follows from Lemma~1.11 in~\cite{VKDiagrams} that, up to homotopy, every map arises in such a way.

\begin{exam}\label{example - pentagon}
     Let $\Gamma$ be the oriented pentagon shown in Figure~\ref{fig - pentagon VK diagram}. By Proposition~\ref{prop - rank 2 implies n partite}, the support of any class representable by a torus in the pentagon is a line consisting of at most two edges. So any second homology class with full support can be represented with no less than three disjoint tori. Let $\al \in H_2(\Sal_\G)$ be the class pictorially shown in the Figure~\ref{fig - pentagon VK diagram}, so~$\supp(\al)$ is the full pentagon. We exhibit a genus two representative for this class using the Van Kampen diagram on the right of the figure. The resulting surface has three vertices (corresponding to the different symbols in the picture), 5 squares and, since every edge gets glued to another edge, it has 10 edges. So the Euler characteristic is $-2$ and hence this diagram indeed represents a connected closed orientable surface of genus two. This gives a map~$f\colon \Sigma_2 \to \Sal_\G$ and the image of a generator of $H_2(\Sigma_2)$ coincides with~$\al \in H_2(\Sal_\G)$---to see this, note the number and orientation of each of the $2$-cells in the Van Kampen diagram.
    \begin{figure}[t]
        \def\svgwidth{0.9\textwidth}
        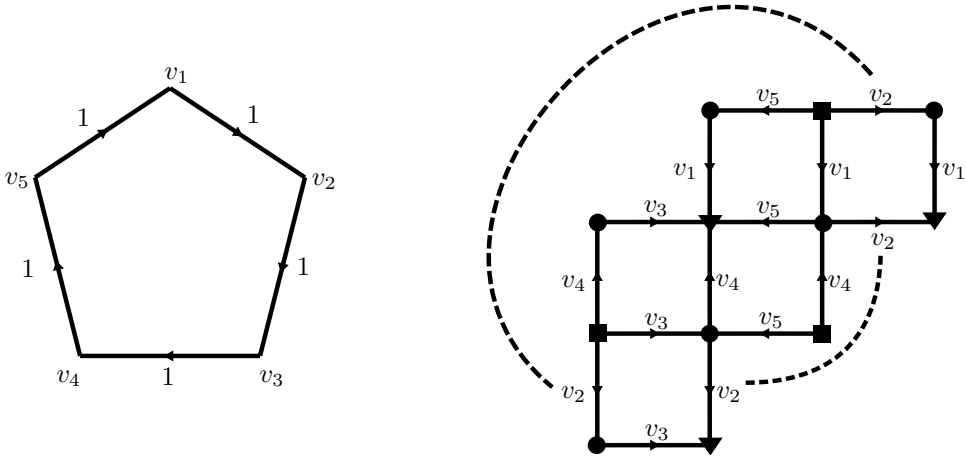
        \caption{On the left is the orientation graph ~$\G$ with labels corresponding to the class $\al \in H_2(\Sal_\G)$. On the right is the Van Kampen diagram corresponding to the map~$f\colon \Sigma_2\to \Sal_\G$, where non-glued opposite edges are identified except for those at the endpoints of the dotted lines which are glued together. Vertices with the same decoration are identified.}
        \label{fig - pentagon VK diagram}
  \end{figure}
\end{exam}
This example provides a negative answer to Question~\ref{question - realised by dj tori}, and thus proves our final Theorem.
\begin{abcthm}\label{abcthm - counterexample to q2}
There exists a RAAG~$A_\G$ with a second homology class whose minimal genus representative cannot be realised by a disjoint union of tori.
\end{abcthm}

\bibliography{genusbib}{}
\bibliographystyle{alpha}
\end{document}